\documentclass{article}
\usepackage{times}
\usepackage{amsmath}
\usepackage[retainorgcmds]{IEEEtrantools}
\usepackage{units}
\usepackage[titletoc,toc,title]{appendix}

\begin{document}
 
\title{Model Vectors}
\author{John Prager\\
   		\texttt{john.prager@gmail.com}}  
\date{\today}  
\maketitle

\newcommand{\seqbegM}{$\textbf{($ }
\newcommand{\seqendM}{$)}$ }
\newcommand{\seqbeg}{\textbf{( }
\newcommand{\seqend}{)} }

\newcommand{\issubsumedby}{\,is \,subsumed \,by \,}

\numberwithin{equation}{section}

\def\changemargin#1#2{\list{}{\rightmargin#2\leftmargin#1}\item[]}
\let\endchangemargin=\endlist

\newcommand{\mysplit}[1]{%
  \begin{tabular}[t]{@{}p{7.5cm}@{}}   
    #1
  \end{tabular}
  }
\renewcommand{\arraystretch}{1.2}

\begin{abstract}
In this article, we discuss a novel approach to solving number sequence problems, in which sequences of numbers following unstated rules are given, and missing terms are to be inferred.  We develop a methodology of decomposing test sequences into linear combinations of known base sequences, and using the decomposition weights to predict the missing term.  We show that if assumptions are made ahead of time of the expected base sequences, then a \emph{model vector} can be created, where a dot-product with the input will produce the result.  This is surprising since it means sequence problems can be solved with no knowledge of the hidden rule.  Model vectors can be created either by matrix inversion or by a novel \emph{combination function} applied to primitive vectors.  A heuristic algorithm to compute the most likely model vector from the input is described.  Finally we evaluate the algorithm on a suite of number sequence problem tests.
  \end{abstract}

\section{Introduction}

\fbox{%
  \parbox{\textwidth}{%
\vspace{1mm}
Arthur, Bob and Claire are trying out party tricks.  Arthur writes on a whiteboard 
\begin{quote}
\begin{tabular}{l l l l l l l l r}
	-2&5&-2&-4&4&&&&?  \\
	\_&\_&\_&\_&\_&&&&\_
\end{tabular}
\end{quote}
and asks Bob to think of a number sequence problem, write 5 terms of the sequence in the blanks under his 5 numbers, and the next term, the missing answer term, under the question-mark.  

\vspace{3mm}

Bob, not one to make it easy for others, is aware that the sequence of squares (1, 4, 9, 16, 25, 36 ...) starts off looking quite similar to the powers of two (2, 4, 8, 16, 32, 64, ...) and supposes that the beginning of the sequence based on their difference, $2^i - i^2$, would look a little mysterious, so he writes these terms.  The whiteboard now shows:
\begin{quote}
\begin{tabular}{l l l l l l l l r}
	-2&5&-2&-4&4&&&&?  \\
	1&0&-1&0&7&&&&28
\end{tabular}
\end{quote}

  }%
}

\fbox{%
  \parbox{\textwidth}{%
\vspace{1mm}
Arthur now asks Bob to multiply together each pair in the 5 left-hand columns, and sum the 5 products.  Bob obeys and computes $-2 + 0 + 2 + 0 + 28 = 28$, which is his next term.  Arthur did not know how Bob would come up with his sequence, and wrote his numbers ahead of time.  ``It's a Magic Multiplier" said Arthur, triumphantly.

\vspace{3mm}

Bob asks ``What if I tried the same sequence, but further to the right?"  ``Try it", says Arthur.  Bob begins with the -1 (which is $2^3 - 3^2$) and enters the next few terms:
\begin{quote}
\begin{tabular}{l l l l l l l l r}
	-2&5&-2&-4&4&&&&?  \\
	-1&0&7&28&79&&&&192
\end{tabular}
\end{quote}
The answer term 192 is  $2^8 - 8^2$, and is correctly predicted.  ``Wow", says Bob.

With Arthur's same 5 numbers, Bob tries again with a decreasing linear sequence 
\begin{quote}
\begin{tabular}{l l l l l l l l r}
	-2&5&-2&-4&4&&&&?  \\
	19&16&13&10&7&&&&4
\end{tabular}
\end{quote}
then an alternating one
\begin{quote}
\begin{tabular}{l l l l l l l l r}
	-2&5&-2&-4&4&&&&?  \\
	10&-3&10&-3&10&&&&-3
\end{tabular}
\end{quote}
both of which are solved with the Magic Multiplier, but then Bob tries a sequence where each term is the sum of the previous two, plus two, starting arbitrarily with (1,2):
\begin{quote}
\begin{tabular}{l l l l l l l l r}
	-2&5&-2&-4&4&&&&?  \\
	1&2&5&9&16&&&&27
\end{tabular}
\end{quote}
but the Magic Multiplier predicts 26, not 27.  ``Oh", says Arthur, ``you have the wrong Magic Multiplier for that:  you want $(-1, 2, 1, -5, 4)$.  ``That works", says Bob, ``and it also works with my linear sequence, but not the others."  At this point Claire steps up to the whiteboard and writes
\begin{quote}
\begin{tabular}{l l l l l l l l l l r}
	2&-3&-5&11&-2&-7&5&&&&?  \\
\end{tabular}
\end{quote}
``That's a Magic Multiplier that works with all your sequences, and then some", she says.  ``But powerful magic needs a powerful energy source; in this case you have to give it 7 terms of your problem sequence." 

\vspace{2mm}
What's going on here?
  }%
}
\vspace{5mm}

One of the usual question types in aptitude tests (or so-called “IQ tests”) is one in which the test-taker must determine the next term in a sequence.  Maybe the most common subtype of these is where the items in the sequence are numbers.  The following sequence segments (where $x$ stands for the desired answer) are simple exemplars of this kind:
\begin{IEEEeqnarray*}{L} 
	 \textbf{... 1,  4,  7,  10,  13,  x ...}  \\
	 \textbf{... 1,  -3,  1,  -3,  1,  x ...}  \\
	 \textbf{... 2,   3,   5,   9,   17,   x ...}  
\end{IEEEeqnarray*}
We use here ellipses as bracketing symbols to indicate that the sequence segment shown is part of a larger (for our purposes, infinite) sequence.

Showing competence at such problems is believed my some educators to indicate a level of intelligence, since the solution is thought to require cognitive skills to discover the underlying rule.  

Associated with each sequence segment is a model (or informally, an explanation).  We are interested here in an automatic solution to such problems: not just finding the value for $x$, but the underlying rule or generating function too.  We will sometimes call this process "recognizing the sequence".
In this paper we will present a novel mathematical approach to recognizing sequences.  We are not attempting to recognize every possible sequence, or indeed every possible sequence that might appear in such tests.  Rather we are presenting a formal analysis that in conjunction with an iterative trial-and-error algorithm can process most test sequence segments in a remarkably simple fashion.
We will define which classes of sequences our approach will handle.  This definition constrains the generating function somewhat, but in practice these restrictions don't appreciably limit our success rate on actual tests.  In this analysis, we only look at problems where the missing term is at the end of the sequence.  While computing missing ``inside'' terms can be marginally trickier for humans, the mathematical problem of rule-fitting is exactly the same, and that is what we are interested in here.

This article is in three parts.  In the first part we describe the theory of \emph{model vectors} (Arthur and Claire's Magic Multipliers), and show two ways by which they can be constructed.  In the second part we dig a bit deeper into the nature of model vectors, and look at heuristics for automatically generating descriptions for human consumption.  In the last part we describe a simple algorithm for applying model vectors to unseen problems, and evaluate it from downloaded aptitude tests.  While an evaluation section may be unusual in mathematical papers, it is the \textit{sine qua non} for experimental sciences.  We are not evaluating the correctness of the mathematical derivations, but the accuracy of the assertion that the models found to underlie test sequences correspond to those of their human creators, and give the same answers that human test-takers are expected to find.

\part{}
\section*{The Mathematical Theory of Model Vectors}
\section{Defining the Problem}
For simpler sequences, the model or explanation behind it can be expressed in just a few words (e.g. ``squares'', ``powers-of-2 plus one'').  The complexity of the model can be used to guide the evaluation of whether it is correct, since all problem sequences have in principle an infinite number of ``next terms". 

We will declare out-of-scope those problems where the number sequence is an incidental attribute of some other sequence, for example the number of letters in the names of the integers, or the number of moons of the planets counting outwards from Mercury, etc.  Also out-of-scope are non-algebraic sequences such as functions of the digits of $\pi$ or the prime numbers.  
We are looking at sequences that have the following properties (and indeed, most test sequences we have seen do have these properties):

Sequence Criteria (SC):
\begin{enumerate}
 \item The $n$th term can be described either as a continuous function of $n$ and/or of $m$ of the previous terms\footnote{By convention, term 1 is the term after the left-hand ellipsis.  $n$ can be any integer.}
 \item While only a finite number of terms are given in any test, the sequence segment can in principle be continued to the right indefinitely.
 \item The sequence segment can be continued backwards indefinitely to the left too, even if that is not commonly done.  So for example, the segment \textbf{... 1,  2,  4,  8 ...} qualifies since it has the “presequence”  \textbf{...  \nicefrac{1}{8},       \nicefrac{1}{4},       \nicefrac{1}{2} ...} .
\end{enumerate}

Typically, sequence terms in observed test questions are integers, although this seems to be entirely for cosmetic reasons, since there is nothing in explanations of real-world test answers (or in our SC) that requires them to be so.  Thus the sequence segment
\begin{quote}
	\textbf{... \nicefrac{1}{2},  \nicefrac{1}{2},  1,  1\nicefrac{1}{2},  2\nicefrac{1}{2},  4 ...}
\end{quote}

passes the SC this way:
\begin{enumerate} \item it uses the Fibonacci rule (each term is the sum of the previous two terms). 
 \item it is clearly infinite to the right.
 \item presequence terms can be simply calculated iteratively from the recurrence relation, so it is infinite to the left too. To illustrate, if we extend 5 places to the left, the sequence looks like
\textbf{... -1\nicefrac{1}{2},  1,  -\nicefrac{1}{2},  \nicefrac{1}{2},  0,  \nicefrac{1}{2},  \nicefrac{1}{2},  1,  1\nicefrac{1}{2},  2\nicefrac{1}{2},  4 ...} 
 \item The formula 
$x_n = x_{n-1} + x_{n-2}$
       is a continuous function of $x_{n-1}$ and $x_{n-2}$
\end{enumerate}

The requirement that the generating function be continuous will appear to eliminate from consideration sequences that use different generating functions depending on the parity (or more generally, any property) of n.  These include sequences like 
\begin{quote}
\textbf{... 2,  3,  6,  7,  14,  15,  30,  31 ...}
\end{quote}
where 
\begin{equation*}
	x_{n+1} = \left\{
	\begin{array}{ll}
	 x_n  +  1  &  \text{if $n$ is odd}, \\
	2x_n &      \text{if $n$ is even}
	\end{array} \right.
\end{equation*}
However, as we will see in Section \ref{alternating}, given the possibility of meta-rules such as  ``functions F1 and F2 alternate'', we can apply the main approach of this paper to the appropriately extracted subsequences to determine the ``inner'' functions.

The previous discussions have revealed a potential nomenclature problem since the word ``sequence'' can be applied either to the numbers provided in a particular test instance, e.g.
\begin{equation}
	\textbf{... 4,  5,  6,  7,  8 ...} \label{integers}
\end{equation}
or to the entire doubly-infinite sequence in which this is embedded, and to which the underlying generating function is equivalent (in the case of \eqref{integers}, the integers).

In what follows, we will call the infinite object the \emph{sequence}, and, when we need to make the distinction, any finite subsection of it such as given in a test a \emph{sequence segment}.  Sequence segments can then be defined by a function, a starting point (or sometimes a subsegment of the sequence) and a length.

\section{Approach}
The approach we will take to recognize sequences automatically will be developed in two stages.  First we set the problem up as one of solving simultaneous linear equations on infinite vectors\footnote{ To make addition of infinite vectors well-defined, we need merely to label one element of each infinite vector as element 1, so that adding and subtracting proceeeds point-wise between the corresponding elements of the participating vectors.  In the examples given here, unless otherwise stated, the first given term of an infinite series after the left-hand ellipsis can be taken as term 1.  The (unshown) terms to the left of it have zero or negative indices.}, and use standard linear algebra to solve them\cite{Strang}.  

We discover that this results in creating a finite vector, which we call a \emph{model vector}, 
and that performing a right-aligned dot-product of the two generates the next term in the sequence, just as Arthur showed.

We then find that there is a novel shortcut to matrix inversion that can be used to find the model vector.

We will use as a motivating example in this section the problem sequence segment
\begin{equation}
	\mathcal{P} = \textbf{... 1,  0,  5,  8,  17,  x ...}  \label{problem}
\end{equation}
This happens to have 5 given elements - this is not critical, but is common in aptitude tests, and for expositional purposes is a nice balance between too few (not enough terms to uniquely identify a simple explanation) and too many (clutter), so we will use this size for our worked examples.  In addition, some tests place the $x$ in the middle of a sequence, but this is not an essentially different problem, so for simplicity of exposition we will concentrate on those where the $x$ is at the right-hand end.

People tend to solve problems like this either by spotting that the sequence is equal to or close to a known sequence, or by differencing consecutive terms.  We might approach this sequence by observing that the numerical trend is not dissimilar to 
\begin{equation}
	\textbf{... 1,  2,  4,  8,  16 ...}
\end{equation}
and so we will subtract that sequence from the former to get:
\begin{equation}
	\textbf{... 0,  -2,  1,  0,  1 ...}
\end{equation}
which does not look terribly promising.  However, further inspection shows that the original sequence also resembles
\begin{equation}
	\textbf{... 0,  1,  4,  9,  16 ...} \label{wrong}
\end{equation}
and doing the subtraction here results in
\begin{equation}
	\textbf{... 1,  -1,  1,  -1,  1 ...} \label{alternating1}
\end{equation}
Using the intuitive notion that “simple is probably correct”, the hypothesis is that the problem sequence here is the sum of two \emph{base} sequences, the squares and constants of alternating sign, and so the next term is found by summing the next terms of both of the base sequences, \textbf{25} and \textbf{-1}, to arrive at \textbf{24}.

\subsection{Sequences and Vectors}
In the rest of this article, we will be treating sequence segments much like vectors: we will be multiplying them by scalars, adding/subtracting them and forming dot-products with actual vectors.  Formally, to do at least the last of these, we need to take our sequence segment and construct from it a vector with the same elements in the same order, and then proceed with the operation with the just-formed vector.  Since making this step explicit would add nothing useful to the exposition, we will not make any distinction between sequence segments and their corresponding vectors; both to indicate this dual nature, and to distinguish them visually from equation numbers, they are displayed in boldface.  When the subsequence is treated as a vector for purposes of the dot-product operation, it will be shown still with ellipses and enclosed by parentheses.

\subsection{Model Vectors}
\label{modelvectors}
We will show now that by making some assumptions or hypotheses about the possible constitution of problem sequences (in this case, of length 5) one can a priori construct a \emph{model vector} such as 
\begin{equation}
	\bar{\mathbf{v}} = \textbf{(-2,  5,  -2,  -4,  4)} \label{mv}
\end{equation}
which has the property that when entered into a dot-product with the given problem segment (and/or many others) the result  is the desired answer.  For example, all of the problem sequence segments in \eqref{works} (and many more) can be solved (i.e. extended one more term) by forming the dot-product with $V$:
\begin{IEEEeqnarray*}{LL}
	\bar{\mathbf{v}}.\textbf{(... 1, 4, 7, 10, 13 ...)} & = \textbf{16} \\   
	\bar{\mathbf{v}}.\textbf{(... 8, 16, 32, 64, 128 ...)}  & = \textbf{256}  \\
	\bar{\mathbf{v}}.\textbf{(... 2, 8, 18, 32, 50 ...)}  & = \textbf{72}  \IEEEyesnumber \label{works} \\
	\bar{\mathbf{v}}.\textbf{(... 6, 5, 6, 5, 6 ...)}  & = \textbf{5} \\
	\bar{\mathbf{v}}.\textbf{(... 1024, 2049, 4098, 8195, 16388 ...)}  & = \textbf{32773} \\
	\bar{\mathbf{v}}.\textbf{(... -9, -4, -1, 0, -1 ...)} &  = \textbf{-4}
\end{IEEEeqnarray*}
We will henceforth refer to this $\bar{\mathbf{v}}$ as our first Example Model Vector, or \textbf{EMV1}.

The observation here is that the problem sequence is a linear sum of other sequences, which we call \emph{base sequences}, that in some sense are “primitive”.  The set of base sequences under consideration at any point we call a \emph{theory}.  For any such theory, we can construct via linear algebra a model vector that embodies this theory. 

\subsection{Matrix Inversion approach}
Let us suppose that our test sequence is constructed from some combination of the base sequences in a set $S$.  For our worked example we will choose the following five 5-element segments, with descriptive labels for easy reference:
\begin{IEEEeqnarray*}{Llll}
	S_1: & \textbf{... 1, 1, 1, 1, 1 ...} & & \text{constant} \\   
	S_2: & \textbf{... 1, 2, 3, 4, 5 ...} & & \text{linear} \\
	S_3: & \textbf{... 1, 4, 9, 16, 25 ...} &\ & \text{quadratic} \IEEEyesnumber \label{s} \\
	S_4: & \textbf{... 1, 2, 4, 8, 16 ...} & & \text{powers-of-2} \\
	S_5: & \textbf{... 1, -1, 1, -1, 1 ...} & & \text{alternating}
\end{IEEEeqnarray*}

There is nothing required or inevitable about these sequences.  They have been chosen here because by observation many test sequences can be analyzed as combinations of these.  In a later section we will address the problem of how to choose the best set of base sequences for a given test sequence.

We call the set of descriptive labels the \emph{description} of a theory, and we will sometimes refer to a theory through its description.  This is discussed further in Section \ref{explaining}.

Because a test sequence might happen to be equal to any of the $S_i$, the $S_i$ must necessarily obey (at least) the same requirements as the test sequence, namely that they obey the previously-defined SC.

The reader might ask, what if the test sequence is indeed composed of the some of the given base sequences, but \emph{shifted} (i.e. taken from some other subsection of their corresponding infinite sequences) - for example \textbf{... 9, 16, 25, 36, 49 ...}, which also can be thought of as representing the \emph{quadratic} sequence?  Since we are analyzing the test sequence as linear combinations of the base sequences, it is sufficient if each base sequence can be shifted (in either direction) through linear combinations of it and other sequences in $S$.  Thus for membership in $S$, not only must the $S_i$ individually satisfy the SC, but they must be shiftable in either direction by linear combination with (some subset of) the other $S_j$.  This means that whether a sequence $S_*$ can be a member of a theory depends on both its own properties and the other sequences in the theory.

\subsubsection{Base Sequence Set Criteria}
\label{bssc}
We therefore introduce the \emph{Base Sequence Set Criteria} (BSSC), to govern which sequences can be used together.
\begin{enumerate}
 \item Every $S_i$ in $S$ satisfies SC
 \item The sequences in $S$ are linearly independent
 \item Any sequence in the set can be shifted any number of positions by forming linear combinations of the sequences in $S$
\end{enumerate}
The second set criterion is required so that any solution is unique.  The third is so that we are not beholden to the particular subsequence of the corresponding infinite sequence that the $S_i$ consist of.   Let us look at each of the $S_i$ in \eqref{s} in turn, to establish what must be done to shift it $k$ elements (for any integer $k$):
\begin{enumerate}
 \item constant.  All subsequences of the infinite constant sequence look the same, so nothing is required
 \item linear.  Add $k$ times $S_1$ to $S_2$
 \item quadratic. Since $(j+k)^2 = j^2 + 2jk + k^2$ add  $2k.S_2$ and $k^2.S_1$ to $S_3$.  This generalizes for all polynomial sequences.
 \item powers-of-2.  Multiply $S_4$ by $2^k$
 \item alternating. Multiply $S_5$ by $(-1)^k$
\end{enumerate}

We suppose that we have a set $S$ of $n$ sequence segments $\{S_i\}$ each of length $n$ elements, and that each of these passes the BSSC set criteria just introduced. We first find a general formula for solving sequence problems by matrix inversion, then apply it to our worked example.

\subsubsection{Creating the Model Vector by Matrix Inversion}
\label{creatingmv}
We now proceed to construct a model vector.
If we assume that the $n$-element problem sequence $\mathcal{P}$ is composed of a linear combination of the $S_i$, then there exists a weight vector $\mathbf{w}$ of length $n$ such that
\begin{equation*}
	 \mathcal{P} = \sum_{i=1}^n w_i S_i
\end{equation*}
Constructing a matrix $M$ whose columns are the $S_i$, we have
\begin{equation*}
	\mathcal{P} = M\mathbf{w}
\end{equation*}
so\footnote{The columns of $M$ are independent due to BSSC-2, so $M$ can be inverted.}
\begin{equation}
	\mathbf{w} = M^{-1}\mathcal{P} \label{inverse}
\end{equation}
We now construct a vector $\mathbf{z}$ of the $n$ ``next'' terms of the $S_i$, where $z_i$ is the next term of $S_i$.  We now can solve the problem since the
solution, i.e. the next term $p_{n+1}$ of the problem sequence $\mathcal{P}$, is 
\begin{equation*}
	p_{n+1} = \mathbf{z}.\mathbf{w}
\end{equation*}
Combining with \eqref{inverse}, we get through associativity
\begin{equation*}
	p_{n+1} = \mathbf{z}.(M^{-1}\mathcal{P}) = (\mathbf{z}. M^{-1}).\mathcal{P}
\end{equation*}
We note that since $M^{-1}$ and $\mathbf{z}$ are independent of $\mathcal{P}$, they can be multiplied in advance of being given $\mathcal{P}$.  Thus, setting
\begin{equation}
	(\mathbf{z}. M^{-1}) = \bar{\mathbf{v}} \label{v} \\
\end{equation}
we get
\begin{equation}
	p_{n+1} = \bar{\mathbf{v}}.\mathcal{P}  \label{x}
\end{equation}
This is the main result, and is the basis for the model vector technique.

When we form the dot-product of $\bar{\mathbf{v}}$ with our problem sequence $\mathcal{P}$, we get the desired next term.  Note that the derivation of \eqref{x} was independent of the specific values of the base sequence segments $S_i$, and in fact of their length.

\subsubsection{Worked Example}
For the 5 sequences in \eqref{s}, the corresponding matrix $M$ is 
\begin{equation*}
\begin{pmatrix}
1 & 1 & 1 & 1 & 1 \\
1 & 2 & 4 & 2 & -1 \\
1 & 3 & 9 & 4 & 1 \\
1 & 4 & 16 & 8 & -1 \\
1 & 5 & 25 & 16 & 1 
\end{pmatrix}
\end{equation*}
[Note that the columns can be arranged in any order, as long as the "next-terms" in \eqref{n} are ordered the same way.]

For the base sequences $S_i$ in \eqref{s}, the next-term vector $\mathbf{z}$ is
\begin{equation}
	\mathbf{z} = \textbf{(1, 6, 36, 32, -1)} \label{n}
\end{equation}
Applying \eqref{v}, we get 
\begin{equation}
	(\mathbf{z}. M^{-1}) = \bar{\mathbf{v}} = \textbf{(-2, 5, -2, -4, 4)} \label{vv}
\end{equation}
This is our \textbf{EMV1}.

When we form the dot-product of $\bar{\mathbf{v}}$ in \eqref{vv} with our problem sequence vector $\mathcal{P}$ from \eqref{problem}, 
\begin{equation}
	\mathcal{P} = \textbf{(... 1,  0,  5,  8,  17 ...) }  \label{problemV}
\end{equation}
we get \textbf{24}, the desired next term.

\subsubsection{Inversion Details}
By using standard matrix inversion methods, we find that 
\begin{equation*}
M^{-1} =
\begin{pmatrix}
11/4 & -13/8 & -13/8 & 17/8 & -5/8 \\
-2 & 2 & 3/2 & -2 & 1/2 \\
1/2 & -3/4 & -1/4 & 3/4 & -1/4 \\
-1/3 & 2/3 & 0 & -2/3 & 1/3 \\
1/12 & -7/24 & 3/8 & -5/24 & 1/24 
\end{pmatrix}
\end{equation*}
with which the reader can verify that \ref{vv} holds.

The weight vector $\mathbf{w}$, which we did not need to find explicitly, is 
\begin{equation}
	\mathbf{w} = M^{-1}\mathcal{P} = \textbf{(1, -2, 1, 0, 1)} \label{w}
\end{equation}
The first three terms have the effect of constructing \eqref{wrong} (the quadratic sequence segment starting from 0), from $S_3$ (the quadratic sequence segment starting at 1), by shifting it one position to the left, according to 
\begin{equation*}
	(m-1)^2 = m^2 -2m + 1
\end{equation*}
giving the segment \eqref{wrong}.  The last term of $\mathbf{w}$ adds the alternating-sign segment \eqref{alternating1}.  We discuss this process in more detail in Section \ref{sec:base}.

\subsubsection{Limitations of the Approach}
We might note that \emph{inter alia} the sequence 
\begin{equation}
	\textbf{... 4, 7, 11, 18, 29 ...} \label{fib}
\end{equation}
(each term is the sum of the previous two) is not correctly solved by \textbf{EMV1} since there is no way to construct the Fibonacci sequence from linear combinations of the sequences in \eqref{s}. That is, the theory embodied by \textbf{EMV1} is not the right one for this test sequence.  This is the same reason that \textbf{EMV1} failed for Bob's fifth challenge, namely the problem sequence
\begin{equation}
	\textbf{... 1, 2, 5, 9, 16 ...}
\end{equation}

Using a different theory, for example one with the following set of sequences $S'$ =
\begin{IEEEeqnarray}{Lll}
	S_1' = S_1: & \textbf{... 1, 1, 1, 1, 1 ...} & \text{constant} \label{sf} \nonumber \\   
	S_2' = S_2: & \textbf{... 1, 2, 3,4, 5 ...} & \text{linear} \nonumber \\
	S_3' = S_3: & \textbf{... 1, 4, 9, 16, 25 ...} & \text{quadratic} \\
	S_4': & \textbf{... 0, 1, 1, 2, 3 ...} & \text{Fibonacci}  \nonumber \\
	S_5': & \textbf{... 1, 1, 2, 3, 5 ...} & \text{Fibonacci (shifted)} \nonumber
\end{IEEEeqnarray}

we can construct a model vector
\begin{equation}
	\textbf{(-1,  2,  1,  -5,  4)}  \label{otherMV}
\end{equation}
which will give the correct solution, namely \textbf{27}.  By adding $S_4'$ and $S_5'$ we get the Fibonnaci sequence shifted one element further, and it is easy to show that adding appropriate multiples of $S_4'$ and $S_5'$ will result in any desired subsequence of the Fibonacci sequence.  Thus \eqref{otherMV} can be used to solve any problem sequence that can be decomposed into a Fibonacci sequence and polynomials of degree up to 2.  We will henceforth refer to it as \textbf{EMV2}.

We note that we had to drop the alternating-sign and powers-of-two base sequences from $S$ in order to keep its size at 5, and hence lose the ability to solve problems based on them.  If we had wanted to keep these two, then the derived model vector $\bar{\mathbf{v}}$ would have 7 elements, as Claire demonstrated, and would require a problem sequence of size at least 7;  in practice we have found most test sequences are not this long.   Our example model vectors \textbf{EMV1} and \textbf{EMV2} can both solve a number of problem sequences, but, although there is overlap, not the same ones.  The limitations of any one particular model vector would seem to raise the question of which one to choose when given a problem sequence.   In Part \ref{partEval} we present an algorithm in which we explore the space of model vectors in order to dynamically select, in fact to create, the one to apply.

\subsubsection{Model vector length}
We have shown the application of a model vector of the same length as the problem sequence, since this is what the mechanics of the dot product requires.  However, if we have a problem that is longer than a model vector, we can simply use the later terms of the problem sequence, ignoring the earlier terms.  This allows the computation to proceed, but is potentially unsatisfactory since it is not taking advantage of information in the earlier terms.  However, if the model vector is indeed appropriate to the problem in hand, then applying it to the leftmost terms of the problem sequence should construct the later terms.  This is the basis of our notion of \emph{consistency} (Section \ref{passing}) and the dynamic solution in Section \ref{sec:dymvec}.

\subsection{Base Sequences}
\label{sec:base}.
We now recap what constitutes a base sequence.

A model vector $\bar{\mathbf{v}}$ is an ordered collection of $n$ real numbers\footnote{Actually, complex numbers, but no non-real numbers, in fact no irrational numbers have been sighted by us in the wild in these number sequence tests} $v_i$.  We will call $n$ its \emph{length}; we will form scalar products with it and problem sequence segments of the same length.  When displayed horizontally we will call the left-most term the \emph{trailing} term and the rightmost term the \emph{leading} term.  Clearly, extending it to the left with zeros will not affect its operation in scalar products; this will be a useful construction later when we need to consider together model vectors of the same length.

We saw above how the problem sequence was deconstructed into a linear combination of the base sequences in \eqref{s} (although the multiples did not have to be explicitly realized in order to compute the next term).  This process would not be very useful if it depended on exactly which sub-sequence of the infinite sequences were used in the construction of $M$.  That is, we would like it to be sufficient to include one representative segment of these infinite sequences, and let the mechanism of the solution process take care of \emph{shifting} the given segment to the appropriate location.

Since our model vector construction is based on linear combinations of base sequences, we automatically deal with problem sequences that are, or are based on, multiples of base sequences - for example 
\begin{equation*}
	\textbf{(3, 3, 3, 3, 3)} = 3 \times \textbf{(1, 1, 1, 1, 1)}
\end{equation*}
However, there is no reason to expect that a given problem sequence will have been formed directly from the particular base sequence segment that is used to construct $M$, as was done with our problem sequence \eqref{problem}.  For example, the sequence of squares 
\begin{equation*}
	\textbf{... 16, 25, 36, 49, 64 ...}
\end{equation*}
is not present explicitly in \eqref{s}.  However, from the identity
\begin{equation*}
	(n+k)^2 = n^2 + 2kn + k^2
\end{equation*}
we see that a \emph{shifted} quadratic sequence is equivalent to any given quadratic sequence plus multiples of linear and constant sequences - which is automatically handled in the case of \eqref{s} since our theory (the set of sequences under consideration) includes the linear and constant sequences.  Put another way, one cannot shift the linear or quadratic sequences by taking multiples of any one alone, but linear combinations of the set \{constant, linear, quadratic\} can produce any of them shifted by any amount.  In general, if a theory $\mathcal{T}$ is to include sequences based on polynomials of degree $n$, then $\mathcal{T}$ must include sample segments of all degrees from 1 to $n$.

\section{Model Composition}

Suppose we have two theories, $\mathcal{T}_a$ and $\mathcal{T}_b$, with corresponding model vectors $\bar{\mathbf{a}}$ and $\bar{\mathbf{b}}$.  Recall that a theory is a set of base sequences.  If we combine the base sequences in $\mathcal{T}_a$ and $\mathcal{T}_b$ we will have a new theory, $\mathcal{T}_c$, say.  $\mathcal{T}_c$ will have a corresponding model vector $\bar{\mathbf{c}}$.  Since we know that the length of a model vector (denoted by $|\bar{\mathbf{v}}|$) is the same as the number of sequences in its theory, we can conclude that $\bar{\mathbf{c}}$ will be of length $|\bar{\mathbf{a}}| + |\bar{\mathbf{b}}|$ (subject to trimming trailing zeroes, which we discuss below).  We ask if we can derive $\bar{\mathbf{c}}$ from $\bar{\mathbf{a}}$ and $\bar{\mathbf{b}}$.  I.e., is there a composition operator $\otimes$ such that $\bar{\mathbf{c}} = \bar{\mathbf{a}} \otimes \bar{\mathbf{b}}$?

This is practically speaking an important question.  If we have a model vector $\bar{\mathbf{v}}$ of length $n$ and we wish to add a single new base sequence (say) to its theory, we would (apparently) have to now invert an $n+1$ x $n+1$ matrix.  What we are about to show is that we can piggy-back on the $n$ x $n$ inversion inherent in $\bar{\mathbf{v}}$ and perform a fairly simple computation on $\bar{\mathbf{v}}$ to come up with the new model vector.  In fact, by building up model vectors from scratch, we will never have to perform any matrix inversion at all.

\subsection{The Model Composition Formula}
\label{sec:comp}

We first observe that in the context of solving test sequences for the next term (to the right),
for any model vector, we can append any number of zero elements to its left, since they will have no effect on dot-product operations (except requiring that formally, for test sequences, enough terms be present).  

We will define here the model composition formula for two model vectors $\bar{\mathbf{a}}$ and $\bar{\mathbf{b}}$ of the same length $n$, where if they were not originally of the same length, the shorter one (wlog $\bar{\mathbf{b}}$) has been padded with zeroes to the left.

We will assert, then prove, that $\bar{\mathbf{c}}$ as defined above is computed by the following sum
\begin{equation}
	\bar{\mathbf{c}} = {_{0n}\bar{\mathbf{a}}} + {_{0n}\bar{\mathbf{b}}} - \mathbf{f} \label{comp_formula}
\end{equation}
where
\begin{quote}
$_{0n}\bar{\mathbf{a}}$ = $\bar{\mathbf{a}}$ prepended with $n$ zeroes.
\end{quote}
\begin{quote}
$_{0n}\bar{\mathbf{b}}$ = $\bar{\mathbf{b}}$ prepended with $n$ zeroes (in addition to any prior padding).
\end{quote}
\begin{quote}
$\mathbf{f}$, whose $m$'th element $(1< m <=2n-1)$ is \footnote{A convolution of the appropriate subsegments of $\bar{\mathbf{a}}$ and $\bar{\mathbf{b}}$}
\begin{equation*}
	\bar{\mathbf{f_m}} = \sum^{j+k=m+1}a_j b_k
\end{equation*}
where $a_j$ is the $j$th element of $\bar{\mathbf{a}}$ (and simultaneously the $j+n$th element of $_{0n}\bar{\mathbf{a}}$), $b_k$ is the $k$th element of $\bar{\mathbf{b}}$ (and the $k+n$th element of $_{0n}\bar{\mathbf{b}}$), $j$ and $k$ range from 1 to $n$ subject to the constraint, and element $f_{2n}$ is zero.  
\end{quote}
We call \eqref{comp_formula} the \emph{Model Composition Formula}.

It is obvious that each of $_{0n}\bar{\mathbf{a}}$, $_{0n}\bar{\mathbf{b}}$ and $\mathbf{f}$ will have (at least) $l$ trailing zeroes.  By inspecting the construction of $\bar{\mathbf{c}}$ we see that if the trailing terms of $\bar{\mathbf{a}}$ and pre-padded $\bar{\mathbf{b}}$ are non-zero, then the $l+1$th term of $\bar{\mathbf{c}}$ will be non-zero too.  So $\bar{\mathbf{c}}$ will have exactly $l$ trailing zeroes, which we will trim for use in practice.  Thus, for the trimmed $\bar{\mathbf{c}}$ and original pre-padded $\bar{\mathbf{b}}$
\begin{equation*}
	|\bar{\mathbf{c}}| = |\bar{\mathbf{a}}| + |\bar{\mathbf{b}}|
\end{equation*}

\subsection{Visualization of the Model Composition Formula}
We can visualize the pattern of the Model Composition Formula by arranging the participating vectors vertically on the page and examining the construction on a term-by-term basis.  We are showing here the rightmost (leading) four terms at the top of the graphic and the leftmost (trailing) two terms at the bottom.
\vspace{3mm}

\begin{tabular}{llll|l}
	Term& $_{0n}\bar{\mathbf{a}}$ & $_{0n}\bar{\mathbf{b}}$ & $\mathbf{f}$ & $\bar{\mathbf{c}}$  \\
\hline
	$2n$ & $a_n$ & $b_n$ & 0 & $a_n + b_n$ \\   
	$2n-1$ & $a_{n-1}$ & $b_{n-1}$ & $a_n b_n$ & $a_{n-1} + b_{n-1} - a_n b_n$ \\
	$2n-2$ & $a_{n-2}$ & $b_{n-2}$ & $a_{n-1} b_n + a_n b_{n-1}$ & $a_{n-2} + b_{n-2} - a_{n-1} b_n -  a_n b_{n-1}$ \\
	$2n-3$ & $a_{n-3}$ & $b_{n-3}$ & $a_{n-2} b_n + a_{n-1} b_{n-1} + a_n b_{n-2}$ & $a_{n-3} + b_{n-3} - a_{n-2} b_n - a_{n-1} b_{n-1}$ \\
	& & & & \hspace{5em} $- a_n b_{n-2}$ \\
\hline
	2 & 0 & 0 & $a_1 b_2 + a_2 b_1$ & $- a_1 b_2 - a_2 b_1$ \\
	1 & 0 & 0 & $a_1  b_1$ & $-a_1 b_1$
\end{tabular}

\subsubsection{Examples} \label {examples}

The following examples might assist in understanding the composition operations discussed in Section \ref{sec:comp}.

For two vectors with single terms
\begin{equation*}
	(a) \otimes (b) = (-a b, a + b)
\end{equation*}

For two vectors where one is single-termed
\begin{equation*}
	(a_1, ..., a_{n-2}, a_{n-1}, a_n) \otimes (b) \\
= (-a_1 b, ..., a_{n-3} -a_{n-2} b, a_{n-2} -a_{n-1} b, a_{n-1} -a_n b, a_n + b)
\end{equation*}

\subsection{Proof of the Model Composition Formula}
\label{proof}

We need to show that if model vector $\bar{\mathbf{a}}$ recognizes sequence $\mathcal{Y}$ and $\bar{\mathbf{b}}$ recognizes $\mathcal{Z}$ then $\bar{\mathbf{c}} = \bar{\mathbf{a}} \otimes \bar{\mathbf{b}}$ recognizes both $\mathcal{Y}$ and $\mathcal{Z}$.  Since $\otimes$ is symmetric it is sufficient to show that $\bar{\mathbf{c}}$ recognizes $\mathcal{Y}$.

We consider a length-$2n$ subsequence of $\mathcal{Y}$, instantiated as vector $\mathbf{y}$.  Since $\mathcal{Y}$ is recognized by $\bar{\mathbf{a}}$ we have for all $0<=k<=n$
\begin{equation}
\sum_{i=1}^na_iy_{i+k} = y_{n+k+1}   \label{sumay}
\end{equation}

This is applying the length-$n$ model vector $\bar{\mathbf{a}}$ to successive length-$n$ subsequence segments of $\mathbf{y}$.  We interpret $y_{2n+1}$ as the ``next term'' of the sequence (i.e. the term to be found).

From the definition of $\bar{\mathbf{c}}$ in \eqref{comp_formula}, we can represent it as the sum of the following vectors (all of length 2n):
\begin{IEEEeqnarray*}{Lll}
& & (0, 0, 0, ..., 0, a_1, a_2, ..., a_n) \\
+& &(0, 0, 0, ..., 0, b_1, b_2, ..., b_n) \\
-&b_n & (0, 0, ..., 0, a_1, a_2, ..., a_n, 0) \\
-&b_{n-1} & (0, ..., a_1, a_2, ..., a_n, 0, 0) \\
& ... \\
-&b_2 & (0, a_1, a_2, ... a_n, 0,..., 0) \\
-&b_1 & (a_1, a_2, ... a_n, 0,..., 0, 0) 
\end{IEEEeqnarray*}
Then for vector $\mathbf{y} = (y_1, y_2 ... y_{2n})$, the dot-product $\bar{\mathbf{c}}.\mathbf{y}$ will be
\begin{equation*}
\sum_{i=1}^na_iy_{i+n}  +\sum_{i=1}^nb_iy_{i+n} - b_n \sum_{i=1}^na_iy_{i+n-1} - b_{n-1} \sum_{i=1}^na_iy_{i+n-2} ...  - b_1 \sum_{i=1}^na_iy_{i-1}
\end{equation*}
So
\begin{equation}
\bar{\mathbf{c}}.\mathbf{y} = \sum_{i=1}^na_iy_{i+n}  +\sum_{i=1}^nb_iy_{i+n} - \sum_{j=1}^nb_{n+1-j} \sum_{i=1}^na_iy_{i+n-j} \label{eqay}
\end{equation}

By putting $k = n+1-j$ in the last term of \eqref{eqay}, we get
\begin{equation}
\bar{\mathbf{c}}.\mathbf{y} = \sum_{i=1}^na_iy_{i+n}  +\sum_{i=1}^nb_iy_{i+n} - \sum_{k=1}^nb_k \sum_{i=1}^na_iy_{i+k-1} \label{eqa}
\end{equation}
Substituting from \eqref{sumay}, we get
\begin{equation}
\bar{\mathbf{c}}.\mathbf{y} = \sum_{i=1}^na_iy_{i+n}  +\sum_{i=1}^nb_iy_{i+n} - \sum_{k=1}^nb_ky_{n+k} 
\end{equation}
hence
\begin{equation}
\bar{\mathbf{c}}.\mathbf{y}= y_{2n+1}
\end{equation}
So $\bar{\mathbf{c}}$ recognizes $\mathcal{Y}$, Q.E.D.

\part{} \label{partEval}
\section*{Model Vectors: Factors, Descriptions and More Examples}
\section{Base Vectors for Base Sequences}
We have seen that our worked example model vector \textbf{EMV1} is created from the base sequences S \eqref{s}.  We have defined the composition operator  $\otimes$ which is used to create model vectors (and corresponding theories) from simpler (shorter) model vectors (theories).  So we ask whether each of the base sequences $S_i$ has its own model vector.  It turns out that we can associate with each what we call a \emph{base vector}, some of which also qualify as model vectors.

We can define base vectors the following way.  Suppose we have a model vector $\bar{\mathbf{v_S}}$ which corresponds to a set of base sequences $S = \{S_i\}$.  If we add the base sequence $Q$ to $S$ forming set $R$, and the resulting set has model vector $\bar{\mathbf{v_R}}$, then the base vector $ \bar{\mathbf{g_Q}}$ for $Q$ is the vector that satisfies
\begin{equation*}
	\bar{\mathbf{v_S}} \otimes \bar{\mathbf{g_Q}} = \bar{\mathbf{v_R}}
\end{equation*}
If $S$ is empty, then $\bar{\mathbf{g_Q}} = \bar{\mathbf{v_R}}$, since the empty vector is the identity under $\otimes$.

Suppose we have what might be the simplest possible problem sequence, which is the constant sequence 
\begin{equation*}
	\textbf{... p, p, p, p, p ...}
\end{equation*}

Clearly the model vector
\begin{equation*}
	\textbf{(1)} \label{constantmv}
\end{equation*}
is necessary and sufficent to generate the next term via dot-product, and so is the base vector that corresponds to the simple theory \{constant\}.  Similarly, 
\begin{equation*}
	\textbf{(2)}
\end{equation*}
is the base vector for the simple theory \{powers-of-two\}, and in fact
\begin{equation*}
	\textbf{(t)}
\end{equation*} 
is the base vector for the theory \{powers-of-t\}, for any $t$.  Also,
\begin{equation*}
	\textbf{(-1)}
\end{equation*} 
is the base vector for the theory \{alternating\}.  It is easily shown that a singleton theory consisting of any of these alone satisfies the \emph{Basic Sequence Set Criteria} defined earlier, and each is a model vector too.

 But what if we ask what is the base vector for the linear base sequence?

It is clear from elementary algebra that a linear sequence can be extended (solved) by the model vector
\begin{equation*}
	\textbf{(-1, 2)}
\end{equation*}
It is also clear that since there are two degrees of freedom in defining a linear sequence, one cannot expect a shorter model vector.  Now, we notice that this model vector can also solve constant sequences, since they are special cases of linear sequences.  But constant sequences can be solved by the model vector \textbf{(1)}, so we ask if there is a vector $\bar{\mathbf{g}}$ such that 
\begin{equation*}
	\textbf{(1)} \otimes \bar{\mathbf{g}} =  \textbf{(-1, 2)}
\end{equation*}
The solution to this is
\begin{equation*}
	\bar{\mathbf{g}} = \textbf{(1)}
\end{equation*}

Thus \textbf{(1)} is the base vector for linear sequences.
This can be interpreted as the delta required to extend the capabilities of a constant theory to a linear one, in exactly the same way that the linear sequence segment $S_2$ does not by itself satisfy the BSSC, but the set $\{S_1, S_2\}$ does.  

\textbf{(1)} can only act as a base vector for the linear sequence as long as the constant base vector (another \textbf{(1)}) is also part of the same set.  We will denote this contingency by a subscript indicating the number of additional copies required to be operative, thus the base vector for the linear sequence is $\textbf{(1)}_1$.  By convention we will drop the subscript 0.

By similar reasoning, the base vector for quadratic is $\textbf{(1)}_2$.  The model vector for the theory \{constant, linear, quadratic\} is \textbf{(1, -3, 3)}, and it is easy to verify that 
\begin{equation*}
	\textbf{(1)} \otimes \textbf{(-1, 2)} = \textbf{(1, -3, 3)}
\end{equation*}

Thus the theory \{constant, linear, quadratic, powers-of-2, alternating\} is represented by the collection of base vectors \{\textbf{(1)}, $\textbf{(1)}_1$, $\textbf{(1)}_2$, \textbf{(2)}, \textbf{(-1)}\}, and 
\begin{equation}
	\textbf{(1)} \otimes \textbf{(1)} \otimes \textbf{(1)} \otimes \textbf{(2)} \otimes \textbf{(-1)} = \textbf{(-2, 5, -2, -4, 4)} = \textbf{EMV1} \label{mod1}
\end{equation}

This can be confirmed by applying the Model Composition Formula step-by-step as in Section \ref{examples}.

\subsection{Compositions of Little Value}
One might ask, if multiples of the base vector \textbf{(1)} are useful, then what benefit might we get from including other multiples in a theory? 
It can be readily shown that multiples of the base vector for the alternating sequence recognize polynomials with alternating signs.  Thus
\begin{equation*}
	\textbf{(-1)} \otimes \textbf{(-1)} = \textbf{(--1, -2)} \text{   which recognizes   }   \textbf{... 1, -2, 3, -4, 5 ...}
\end{equation*}
and
\begin{equation*}
	\textbf{(-1)} \otimes \textbf{(-1)} \otimes \textbf{(-1)} = \textbf{(--1, -3, -3)} \text{   which recognizes   }   \textbf{... 1, -4, 9, -16, 25 ...}
\end{equation*}

Beyond this, there does not seem much value in using multiples.  Take \textbf{(2)}:  what do we get if we include two of them?  Now, 
\begin{equation*}
	\textbf{(2)} \otimes \textbf{(2)} = \textbf{(-4, 4)} 
\end{equation*}
It is easily confirmed that the model vector \textbf{(-4, 4)} recognizes any sequence that \textbf{(2)} does, but more generally any where each term is four times the difference between the previous two.  This seems to be of very narrow utility in the world of sequence problems.  Likewise, other base vector multiples do not seem particularly useful.  Therefore, when in Section \ref{sec:dymvec} we describe an algorithm for building model vectors up from scratch from smaller base vectors, we might limit the number of multiples of most base vectors to 1.

\subsection{Explaining a Theory}
\label{explaining}
If, loosely speaking, the best theory to fit a sequence is the one with the best (read: shortest) description, it behooves us to develop a notion of a theory description.  We have seen the association of descriptive labels with base sequences/vectors, so a simple description of a theory would be an enumeration of the labels of its component base vectors.  In some cases, model vectors derived by composition have their own descriptors, but not always.  For example 
\begin{equation*}
	\textbf{(1)} \otimes \textbf{(1)} = \textbf{(-1, 2)} 
\end{equation*}
can be given the label ``linear'', but
\begin{equation*}
	\textbf{(2)} \otimes \textbf{(-1)} = \textbf{(2, 1)} 
\end{equation*}
which combines alternating and powers-of-2, doesn't admit readily of a simpler descriptor than just that.  We can develop a lookup table of substitutions such as Table \ref{tab:complabels}, below, which can be used to shorten theory descriptions when the corresponding factors exist.

\begin{table}[hbp]
\begin{tabular}{|l|l|l|}
\hline
	Factors & Composition & Label  \\
\hline
	$\textbf{(1)} \otimes \textbf{(1)}$ & \textbf{(-1, 2)} &  linear \\   
\hline
	$\textbf{(1)} \otimes \textbf{(1)} \otimes \textbf{(1)}$ & \textbf{(1, -3, 3)} & quadratic \\   
\hline
	$\textbf{(1)} \otimes \textbf{(1)} \otimes \textbf{(1)} \otimes \textbf{(1)} $ & \textbf{(-1, 4, -6, 4)} & cubic \\   
\hline
	$\textbf{(1)} \otimes \textbf{(1)} \otimes \textbf{(1)} \otimes \textbf{(1)} \otimes \textbf{(1)}$ & \textbf{(1, -5, 10, -10, 5)} & quartic \\   
\hline
	$\textbf{(1)} \otimes \textbf{(1)} \otimes \textbf{(1)} \otimes \textbf{(1)} \otimes \textbf{(1)} \otimes \textbf{(1)}$ & \textbf{(-1, 6, -15, 20, -15, 6)} & quintic \\   
\hline
	$\textbf{(-1)} \otimes \textbf{(-1)}$ & \textbf{(-1, -2)} & alternating-linear \\   
\hline
\end{tabular}
  \caption{Labels for selected base vector compositions}
  \label{tab:complabels}
\end{table}

The last of these is maybe not standard, but is meant to illustrate that any choice of descriptor will work as long as it is understandable.  [The appearance of binomial coefficients should not be a surprise.]

A theory's description is not the same as an explanation.  The sequence \textbf{... 4, 5, 7, 11, 19, 35 ...} is solved by $\textbf{(-2, 3)} = \textbf{(1)} \otimes \textbf{(2)}$, which has description \{constant, powers-of-2\}, but this is a bit wanting as an explanation.  A reasonable explanation might be that the $n$th term is $3 + 2^n$, or simply ``powers of two, plus 3''.  In other words, an explanation employs not just a list of the base sequences, but what multipliers are used.  These multipliers (weights) can be recovered if desired by the matrix inversion approach.  In the evaluation section later, when a problem is solved by a model vector, we present just the theory's description, not an explanation.

\subsection{Subsumption}
If model vector $\bar{\mathbf{a}}$ is a factor (under model composition) of $\bar{\mathbf{g}}$ we say that $\bar{\mathbf{g}}$ \emph{subsumes} $\bar{\mathbf{a}}$.  This means that any problem sequence $\mathcal{P}$ of length $|\bar{\mathbf{g}}|$ or greater that is solvable by $\bar{\mathbf{a}}$ is also solvable by $\bar{\mathbf{g}}$.  This is easily seen through the matrix-inversion approach.  Since the set of base vectors of $\bar{\mathbf{g}}$ includes the set of base vectors of $\bar{\mathbf{a}}$, and $\mathcal{P}$ is representable by a linear combination of the base vectors of $\bar{\mathbf{a}}$, then it can clearly be equivalently represented by $\bar{\mathbf{g}}$.
Thus, if we look at model vectors for solving polynomials of increasing degree, we see that
\begin{IEEEeqnarray*}{Ll}
	\textbf{(1)} & \issubsumedby \textbf{(-1, 2)} \\
	& \issubsumedby \textbf{(1, -3, 3)} \\
	& \issubsumedby \textbf{(-1, 4, -6, 4)} \\
	& \issubsumedby \textbf{(1, -5, 10, -10, 5)} \\
	& ...
\end{IEEEeqnarray*}
We can use the existence of subsumption to provide an alternative rendering of \eqref{mod1}, namely
\begin{equation}
	\textbf{(1, -3, 3)} \otimes \textbf{(2)} \otimes \textbf{(-1)} = \textbf{(-2, 5, -2, -4, 4)}  \label{mod2}
\end{equation}
which can be interpreted as the theory \{quadratic, powers-of-2, alternating\}, where constant and linear are implicit.

\subsection{Theory size}
When solving number sequence problems and faced with multiple solutions, the simpler one is usually preferred.  A proxy for simplicity is the ``size'' of the corresponding theory.  We saw above how there can be multiple factorizations of a model vector, and hence multiple candidates for theory size.  We also saw how descriptions can be shortened when combinations (compositions) of base vectors have their own descriptive labels. We then assert that the size of a theory is the number of components in its shortest such description.

Thus the size of $\textbf{(-2, 5, -2, -4, 4)}$ is 3, via \eqref{mod2}, not 5 as \eqref{mod1} might suggest.  This apparent shortening reflects the assumption that descriptive labels such as ``quadratic'' are more easily consumable to average solvers than their decomposition might be \footnote{This shortening technique assumes that short labels such as ``quadratic'' are sufficient to describe sequences such as \textbf{(... 1 2 5 10 17 26 37 ...)} (squares-plus-one).}

\subsection{The Fibonacci sequence and other recurrence relations}
A common device used in numerical sequence problems is to have a term equal the sum of the previous two, or some variation of this.  The simplest example of this is probably the Fibonacci sequence
\begin{equation*}
	\textbf{... 0, 1, 1, 2, 3, 5, 8, 13 ...}
\end{equation*}
This suggests using a base sequence such as 
\begin{equation}
	\textbf{... 0, 1, 1, 2, 3 ...} \label{fib1}
\end{equation}
Now, including this in a theory made from sequences such as those we have previously encountered will break the BSSC, since there will be no way to shift \eqref{fib1} with linear combinations of it and the other base sequences.  However, if we instead include both \eqref{fib1} and \eqref{fib2}
\begin{equation}
	\textbf{... 1, 1, 2, 3, 5 ...} \label{fib2}
\end{equation}
(as was done in \eqref{sf}) we will be in good shape, since appropriate linear combinations of this pair will result in any desired sequence subsegment of the Fibonacci sequence.  These two sequence segments do not have 
meaningful base vectors\footnote{Formally, if we solve $a \otimes b = (1,1)$ we find $a = (\phi), b = (1-\phi)$, where $\phi$ is the familiar golden ratio.  It can be demonstrated that adding multiples of their corresponding powers-of sequences generates sequences obeying the Fibonacci rule, but practically speaking it may be most useful to ignore these quark-like factors and just consider their product.}, but taken as a pair they have the base vector/model vector
\begin{equation}
	\textbf{(1, 1)}
\end{equation}
Thus to complete the example in \eqref{sf}, the theory \{constant, linear, quadratic, Fibonacci, Fibonacci shifted\} can be represented by 
\begin{equation}
	\textbf{(1)} \otimes \textbf{(1)} \otimes \textbf{(1)} \otimes \textbf{(1, 1)} = \textbf{(1, -3, 3)} \otimes \textbf{(1, 1)} = \textbf{(-1, 2, 1, -5, 4)}
\end{equation}
which again is our second example model vector, \textbf{EMV2}.

\subsubsection{Variations}
Some simple recurrence relations and corresponding model vectors are shown in Table \ref{tab:examplemvs}.
With each is a sample problem segment it will solve (all beginning with 1 and with simple weights in the matrix-inversion view, although an infinite number of possibilities exist).  

\begin{table}[hbp]
\begin{tabular}{|l|l|p{2cm}|l|p{2cm}|}
\hline
	Model Vector & Factors & Label & Example & Comment  \\
\hline
	\textbf{(1)} & & constant & \textbf{... 1,1,1,1,1 ...} & $p_n = k$\\
\hline
	\textbf{(-1, 2)} & $\textbf{(1)} \otimes \textbf{(1)}$ &  each term is twice the previous minus the one before & \textbf{... 1,2,3,4,5 ...} & = linear.  Subsumes previous \\   
\hline
	\textbf{(1,1)} & &  each term is sum of the previous 2 & \textbf{... 1,1,2,3,5 ...} & $p_n = p_{n-1}+ p_{n-2}$ = Fibonacci \\
\hline
	\textbf{(1,1,1)} & & each term is sum of the previous 3 & \textbf{... 1,1,1,3,5 ...} & $p_n = p_{n-1} + p_{n-2} + p_{n-3}$ \\
\hline
	\textbf{(1,1,1,1)} & & each term is sum of the previous 4 & \textbf{... 1,1,1,1,4 ...} & $p_n = p_{n-1} + p_{n-2} + p_{n-3} + p_{n-4}$\\
\hline
	\textbf{(-1,1)} & & each term is difference of the previous 2 & \textbf{... 1,2,1,-1,-2 ...} & $p_n = p_{n-1} - p_{n-2}$\\
\hline
	\textbf{(1,-1)} & & each term is difference of the previous 2 & \textbf{... 1,2,-1,3,-4 ...} & $p_n = p_{n-2} - p_{n-1}$ \\
\hline
	\textbf{(t)} & & powers-of-t & \textbf{... 1,$t$,$t^2$,$t^3$,$t^4$ ...} & $p_n = t^n$ \\
\hline
\end{tabular}
  \caption{Sample model vectors representing simple recurrence relations}
  \label{tab:examplemvs}
\end{table}

\part{} \label{partEval}
\section*{An Evaluation of the Theory}

\section{Developing the Algorithm}
We describe here how to build a simple program to solve number sequences using the theory of model vectors.  Now, a given sequence segment can be completed in an infinite number of ways, and each completion can be explained in an unlimited number of ways, but intuitively some are better, simpler than others. For our purposes, solving a sequence amounts to determining the theory behind it, represented by a model vector.  A theory can be viewed either as a set of base vectors, or if some are pre-multiplied as in \eqref{mod2}, a set of model vectors.  The act of selecting a theory is the act of selecting its model vector factors.  

There is no such thing as the one ``correct'' answer, since that depends on the author and requires us to depart mathematics for psychology.   We will suppose that a set of base sequences is initially chosen, so for our purposes the best solution for any test sequence is the simplest combination of these base sequences (according to a metric to be described) which is consistent (Section \ref{passing}, below) with the given sequence.

There are basically two approaches to this task.  We could experiment with as large a collection of test data that can be collected, find the model vector that correctly solves the most of these (perhaps \textbf{EMV1} or \textbf{EMV2}), and use that exclusively.  This approach has the simplicity of selecting the tool once and for all, independently of the problem, and can be useful as a parlour trick, as Arthur demonstrated.  Alternatively, a dynamic approach that builds up the smallest possible model vector that seems to explain the given terms in the sequence can perform better, and is described here.

\subsection{Passing and Consistency} \label{passing}
If we have a problem test sequence $\mathcal{P}$ of length $n$, and a model vector of length $n-1$, and when applied to the first $n-1$ terms of $\mathcal{P}$ it predicts the $nth$ term, we say that the model vector \emph{passes}.  Vectors of length $n-i (1 <= i <= n-1)$ can be tested $i$ times (there are $i$ contiguous subsequences of $\mathcal{P}$ of length $n-i$ that do not include the $nth$ term);  if they pass every time they are said to be \emph{consistent} with $\mathcal{P}$.

If we have a test sequence $\mathcal{P}$ of length $n$, then clearly no model vectors of length $n+1$ or greater can be used to predict the next term.  Any vector of length $n$ will produce a result (with inherent theory and explanation), but there is no direct way to get any confidence that the theory is right (one possible indirect way is to derive and try to assess the complexity of the weight vector).  Any shorter vector that is consistent is a candidate for selection.  Successively shorter model vectors must pass more tests in order to be consistent, so the shortest consistent vectors are preferred since there is more evidence that they are correct.

While the preceding statement is an intuitive one, it expresses a belief that in aptitude tests simplicity (of explanation) is a valued feature.  It can be tested (statistically) by gathering large numbers of test questions and confirming that, where alternative solutions exist, the simpler one is most often preferred.  Indeed, in a later section of this paper we deploy the algorithm being developed here against a collection of test sequences, show that it works well, and confirm that the reasons for the small number of failures are not related to this issue.

\subsubsection{Passing Example} \label{passingexample}
A model vector will always predict the next element, but it might not be the one intended by the author.  A good example of that was given in the opening scenario, where Bob's fifth example was based on the Fibonacci sequence, but \textbf{EMV1}'s theory  did not incorporate Fibonacci, so it was doomed.  The first-given answer, \textbf{26}, corresponded to the theory of \textbf{EMV1}, namely the base sequences $S_i$ in \eqref{s}, with weight vector
\begin{equation*}
	\textbf{(1/2, -1/2, 1/2, 1/3, 1/6) } 
\end{equation*}
which corresponds to an explanation probably beyond the top-of-the-head skills of a typical party-goer.

The problem sequence
\begin{equation}
	\textbf{... 1, 2, 5, 9, 16, 27 ...} \label{probseqfib}
\end{equation}
can be constructed from two (different) Fibonacci base sequences plus the constant base sequence.  The simplest model vector to recognize \ref{probseqfib} can therefore be constructed from 
\begin{equation}
	\textbf{(1, 1)} \otimes \textbf{(1)} = \textbf{(-1, 0, 2)}  \label{simplest}
\end{equation}
We will note that \eqref{simplest} is shorter than and is subsumed by Arthur's second model vector, \textbf{EMV2}.  For completeness, 
\begin{equation}
	\textbf{(-1, 0, 2)} \otimes \textbf{(1)} \otimes \textbf{(1)} =  \textbf{(-1, 2, 1, -5, 4)} = \textbf{EMV2}
\end{equation}
\eqref{simplest} not only computes the "next term" \textbf{27} when applied to the previous three terms in \eqref{probseqfib}, but when applied to any three consecutive terms will compute the next one, and so is consistent with it.

\subsection{Deriving the Weight Vector}
\label{weightvector}
We saw in Section \ref{creatingmv} how the weight vector can be found by applying the inverted base sequence matrix to the problem vector.  It is also derivable without matrix inversion, but in a term-by-term fashion, by using model vector factors and consistency tests, which we sketch here.

The idea is that for any of the base sequences $S_i$ in the model vector's theory, if its weight in the weight vector is $w_i$ and we subtract $w_i$ multiples of it from the problem vector $\mathcal{P}$, the resulting problem sequence $\mathcal{P}'$ will be independent of that base sequence, so a model vector $\bar{\mathbf{v}}'$ made from the correspondingly reduced theory will recognize $\mathcal{P}'$.  Moreover, since $|\bar{\mathbf{v}}'|$ will be less than $|\bar{\mathbf{v}}|$, $\bar{\mathbf{v}}'$ will be consistent with $\mathcal{P}'$.

For pairs (or more) of base sequences that are mutually dependent in a theory, such as two instances of the Fibonacci sequence, the same approach is applied except that two (or more) weights are determined together.

In our worked example with model vector \textbf{EMV1} and problem sequence
\begin{equation*}
	\mathcal{P} = \textbf{(... 1,  0,  5,  8,  17 ...) }
\end{equation*}
we found the weight vector was 
\begin{equation*}
	\mathbf{w} = M^{-1}\mathcal{P} = \textbf{(1, -2, 1, 0, 1)} 
\end{equation*}
Let us take the alternating-sign sequence, which has weight 1.  The model vector without this base sequence is 
\begin{equation}
	\textbf{(1)} \otimes \textbf{(1)} \otimes \textbf{(1)} \otimes \textbf{(2)} = \textbf{(1, -2, 0, 2)}  \label{no-alt}
\end{equation}
(loosely, this is $\textbf{EMV1} \div \textbf{(-1)}$, although we haven't formally defined model inverse composition).
So we ask for what value $h$ is \eqref{no-alt} consistent with the result of 
\begin{equation*}
\mathcal{P} - h\textbf{(1, -1, 1, -1, 1)}.
\end{equation*}
It can be easily verified the solution is $h=1$.  This operation can be iterated on successively reduced versions of the original model vector and successively subtracted versions of the problem sequence, to determine all the weights.  

\subsection{Interleaved Sequences}
\label{alternating}
Numerical sequence aptitude tests in practice are not bound by criteria such as our Sequence Criteria or indeed any discernable constraints.  Some take advantage of the visual appearance of the numbers, using properties of the base-10 system, and don't always play by the ``rules'' of reducing fractions to their lowest terms.  For example, the following sequences:
\begin{equation}
	\textbf{... 0.1, 0.2, 0.3, 0.4, 0.5, 0.6, 0.7, 0.8, 0.9, 0.10, 0.11, 0.12, 0.13 ...} \label{decimal_digits}
\end{equation}
and
\begin{equation}
	\textbf{... 1/6, 2/8, 3/10, 4/12, 5/14 ...} \label{fractions}
\end{equation}
are very easy for most people to complete, but would require extensions to our theory to be solved automatically by it.  We call such sequences \emph{visual} since, although each term may be described by a mathematical function of its position in the list, it is the format or rendering on the page that gives the solver the biggest clue.  Solving visual problems is an interesting challenge, but not one we take on here. 

There is one commonly encountered variation, in fact one that is essentially the same as embodied in both  \eqref{decimal_digits} and \eqref{fractions}, that can be readily accommodated, and this is the use of interleaved sequences.  A sequence such as
\begin{equation}
	\textbf{... 1, 100, 2, 200, 3, 300, 4, 400 ...} \label{interleaved}
\end{equation}
requires very different logic when passing from odd to even terms than from even to odd, but when viewed as two interleaved sequences the formula, and hence the model vector(s), for continuing the sequence is much more straightforward to compute.

In practice, the question arises of when to interpret a problem sequence as being of the interleaved type.  This is especially true of relatively short sequences, such as the typical 5-element sequence, since the component 3- and 2-element sequences  usually have individually simple explanations, although taken as a whole are not right.  For example, the quadratic sequence
\begin{equation}
	\textbf{... 1, 4, 9, 16, 25 ...} 
\end{equation}
can be ``explained'' as the two subsequences \textbf{... 1, 9, 25 ...} and \textbf{... 4, 16 ...}, but chances are that in the context of an aptitude test the desired next term is \textbf{36}, not \textbf{28} (if the shorter sequence was \{linear\}) or \textbf{64} (if the shorter sequence was \{powers-of-4\}).  We therefore adopt the simple heuristic of only proposing interleaving if both subsequences are consistent with (in the sense of Section \ref{passing}) the same model vector.

Intuitively, the simpler the connection between the interleaved sequences (there is a simple one in \eqref{interleaved}), the more likely it is intended.  Evaluating such simplicity is beyond the scope of the current paper, so we adopt a simple expedient in the deployed system we evaluate here.  We attempt to solve the problem both considering it to be interleaved and not, and either way may return a prediction.   If only one does, we use that.  Otherwise, if the sum of the sizes of the two interleaved theories is less than the size of the non-interleaved theory, we use the interleaved prediction, else the non-interleaved prediction.

\subsection{DyMVeC Algorithm}
\label{sec:dymvec}
In this section we present the DyMVeC algorithm, Dynamic Model Vector Construction.  In this algorithm, we build up successively longer model vectors until we find one that is consistent with the problem, if possible.  We grant that this might be a bit much to expect any but the most committed solver to execute by hand, but note that each iteration of model vector construction is simpler than the corresponding matrix inversion.

We start with a set $V$ of base vectors $\{\bar{\mathbf{v_i}}\}$.  For each base vector $\bar{\mathbf{v_i}}$ we assign a positive integer $q_i$ which corresponds to the maximum number of multiples of $\bar{\mathbf{v_i}}$ we are interested in using.  For example, we might assign $4$ to the base vector \textbf{(1)}, indicating we are only interested in polynomials up to degree 3.  

Our approach is to find or construct a theory which is consistent with our problem sequence $\mathcal{P}$, and then use it to predict the next term $p_{n+1}$. 

We generate $V_1$ which is the set of all $\bar{\mathbf{v_i}}$ in $V$ of length 1 and test to see if they are consistent with $\mathcal{P}$.  If there is one, then that is chosen.  If there is more than one, then we choose one of the consistent $\bar{\mathbf{v_i}}$ at random (or if the descriptive labels in the theory have a preference ordering, use that).  If there is no consistent element of $V_1$, we then generate $V_2$, which is the set of all $\bar{\mathbf{v_i}}$ in $V$ of length 2, plus all binary combinations (under $\otimes$) of length-1 base vectors.  For any $\bar{\mathbf{v_i}}$, these combinations will include $\bar{\mathbf{v_i}} \otimes \bar{\mathbf{v_i}}$ as long as $q_i >= 2$.

The members of $V_2$ are now tested for consistency as before.  If none is consistent with $\mathcal{P}$, we generate $V_3$ analogously, and then on if necessary until $V_{n-1}$. 

If this process continues to the end without generating a consistent model vector for the entire sequence, we re-interpret it as two interleaved sequences, and see if the same model vector can be found for both.  At the end of this process we will either have generated a consistent model vector of length $1 <= l <= n-1$, or else admit failure.

An alternative to giving up is to try a statistical approach.  Based on experimental trials with large numbers of test sequences, we can determine ``favorite'' model vectors, namely those that empirically are successful with the greatest number of problems.  Then we will use the favorite of the same length as the current test sequence to generate the next term, although no consistency tests will have been passed.  The model vectors \textbf{EMV1} and \textbf{EMV2} are examples of ones that work well in practice.  In an informal study performed with colleagues who were asked to write down a number sequence problem of size 5 from an imagined test, both model vectors correctly predicted the next term in 18 out of 30 cases (but not the same 18).

\section{Evaluation}
For this evaluation we prime DyMVeC with the ``stable'' of base vectors in Table \ref{tab:stable}.  The vectors, maximum multiplicity, and descriptive labels are given.

\begin{table}[hbp]
\begin{tabular}{|l|l|p{4.2cm}|l|}
\hline
	Base Vector & Multiples & Labels & Comment \\
\hline
	\textbf{(1)} & 4 & constant/linear/quadratic/cubic &\\
\hline
	\textbf{(-1)} & 2 & alternating/alternating-linear & \\   
\hline
	\textbf{(1,1)} & 1 &  each term is sum of the previous 2 & Fibonacci \\
\hline
	\textbf{(1,1,1)} & 1 & each term is sum of the previous 3 &  \\
\hline
	\textbf{(1,1,1,1)} & 1 & each term is sum of the previous 4 &  \\
\hline
	\textbf{(-1,1)} & 1 & each term is difference of the previous 2 & \\
\hline
	\textbf{(1,-1)} & 1 & each term is negative difference of the previous 2 & \\
\hline
	$\textbf{(t)}$ & 1 & powers-of-t & for all integer t in 2-10 \\   
\hline
	$\textbf{(1/t)}$ & 1 & powers-of-$1/t$ & for all integer t in 2-10 \\   
\hline
\end{tabular}
  \caption{Set of base vectors, with multiplicity, used in the evaluation}
  \label{tab:stable}
\end{table}

This collection was developed by testing DyMVeC on problem sequences gathered from colleagues and the Web (but not those in the formal test set, described below).
It should be obvious that other choices are possible - for example there is no reason that $t$ should stop at 10 for the powers-of-$t$ base vectors.  Increasing the size of this collection increases the compute time (a minor consideration because for the current configuration the response is near-instantaneous on a standard home computer), but also increases the likelihood of a false-positive result.  Since aptitude tests are set for ordinary people, not necessarily mathematically-inclined, and since more ``elaborate'' base sequences were not encountered in the development stage, the set of selectable bases was left as above. 

We retrieved a number of on-line number sequence tests\footnote{from http:www.fibonicci.com/math/number-sequences-test/ and http://www.queendom.com/tests/quiz/index.htm?idRegTest=756/7} and ran DyMVeC on them.  The results, with some analysis, for each test problem are shown in the Appendix.

\subsection{Discussion}
The evaluation set is unfortunately small, but most numeric sequence aptitude tests available on the Web require purchase, as part of a test-prep program.  However, with the 62 we acquired, our algorithm solved 45 correctly - 73\%.  Remarkably, if armed with just \textbf{EMV1} a test-taker could get 45\%; these tests are multiple-choice, so with 5 offered answers, the rate one would expect to get by simply guesssing would be 20\%, and if not multiple-choice but write-in, random guessing would get close to 0\%.

\subsubsection{Error cases}
The current theory does not handle sequences where there is a non-constant multiplicative factor between successive terms, or some of the more complicated cases of interleaved sequences, or cases with three interleaved sequences.  Problems that were missed that way were \{\#11, \#12, \#25, \#27, \#31, \#57, \#58\}.

In at least one observed case, (\#24), more than one model vector with the same size consistent with the input was generated, and  DyMVeC happened to choose the wrong one. With many more data points we could machine-learn a more sophisticated weighting scheme for model vector factors, to give better-performing estimates of size.

Some problem sequences played upon the representation in base-10.  These are tagged in the tables as VISUAL, and include \{\#13, \#32, \#48, \#54, \#59, \#61, \#62\}.

\#34 and \#40 do not seem to have a satisfyingly simple interpretation.

\subsubsection{Alternative view}
It might be noted that in some sense the DyMVeC algorithm is never completely wrong.  Given a particular number sequence problem, any solution it proposes is accompanied by a model vector with an explanation, namely a set of base sequences with multipliers.  If the model vector is shorter than the problem sequence, then it has been shown to be consistent, (in the sense of Section \ref{passing}).  Thus even though the explanation might not be the same as the setter's, and may well be subjectively less simple, it does at least explain the input and proposes a continuation accordingly.  If the model vector is of the same size as the problem, then a weaker version of the same argument can be made:  there is a breakdown of the problem sequence into linear combinations of simple sequences, and this analysis has the given solution.

\section{Summary}
We have presented here the theory of \emph{model vectors}, a mechanism for solving certain kinds of number sequence problems found in aptitude tests.   Given a model vector, the next term of a sequence can be found by as simple an operation as a dot-product.  As an illustration, we have shown in Section \ref{modelvectors} how a single model vector, our example vector \textbf{EMV1}, can solve a variety of such problem sequences.  

Under the hypothesis that problem sequences are linear combinations of commonly-preferred infinite base sequences, a way is shown to derive model vectors by matrix inversion.  While mathematically straightforward, this is not an easy thing to do by pen-and-paper for model vectors (or problem sequences) of any reasonable length.  We therefore show an alternative approach, whereby very simple base vectors are associated with the infinite base sequences, and a model composition function is defined by which base vectors can be composed to form the model vectors.  

Model vector lengths are equal to the number of base sequences they represent, so the more capable a model vector is, the longer the problem sequence must be to allow a dot-product.  Since a solver has no say over the length of the problem sequence, which in practice is commonly 5 terms, the solver can benefit from being able to generate a ``most likely'' model vector on-the-fly, based on redundancy within the problem sequence.  Doing so requires acceptance of the proposition that the simplest solution is the best, and also some way to measure the complexity (or simplicity) of a model vector.

With a method to develop a sequence of model vectors, simplest first, plus a test of whether a model vector is consistent with a problem, we demonstrate an algorithm to solve number sequence problems with dynamic model vectors, DyMVeC.  We evaluate DyMVeC on unseen sequences downloaded from the Web.  The evaluation is not an evaluation of the correctness of the mathematical theory, but an evaluation of the hypothesis that number sequence problems can be deconstructed as linear combinations of appropriate base sequences, in a way that aligns with users' understanding.  Our algorithm gets a score of 73\% and the static example vector \textbf{EMV1} gets a score of 45\%; by contrast, guessing in a multiple-choice context would get only 20\%.  What is remarkable is that in neither case is there any deep ``understanding'' of the nature of the problem.

\section*{Acknowledgements}
I would like to thank Roger Heath-Brown for feedback on a very early draft, and my colleagues in the Watson Lab at IBM Research for providing sample problem sequences.  I also want to thank Greg Dresden and Scott Gerard for comments on the current version of the paper.

\pagebreak
\begin{appendices}
\section{Evaluation Details}
We present here the results of running DyMVeC on the numerical sequence problems we found on www.fibonicci.com and www.queendom.com.  These were first seen and downloaded \textbf{after} the development of the model vector theory, and the development and debugging of the DyMVeC algorithm. 

There were 62 such problems, of which 45 (73\%) were correctly solved.  The non-dynamic example model vectors  \textbf{EMV1} and \textbf{EMV2} scored 28 and 19 (45\% and 31\%) respectively.  Random guessing in a 5-choice multiple-choice context would score 20\%, on average.  
\begin{table}[hbp]
\hspace*{-3.5cm}
\begin{tabular}{|l|p{4.5cm}|l|l|l|p{7.5cm}|l|l|}
\hline
	& Problem & Ans & Ours & OK & Explanation & EMV1 & EMV2 \\
\hline
1 & \textbf{2, 4, 9, 11, 16} & 18 & 18 & Y & \{linear, alternating\} & Y & N \\
2 & \textbf{30, 28, 25, 21, 16} & 10 & 10 & Y & \{quadratic\} & Y & Y \\
3 & \textbf{-972, 324, -108, 36, -12} & 4 & 4 & Y & \{reciprocal-powers-of-9\}, interleaved sequences & N & N \\
4 & \textbf{16, 22, 34, 52, 76} & 106 & 106 & Y & \{quadratic\} & Y & Y \\
5 & \textbf{123, 135, 148, 160, 173} & 185 & 185 & Y & \{linear, alternating\} & Y & N \\
6 & \textbf{0.3, 0.5, 0.8, 1.2, 1.7} & 2.3 & 2.3 & Y & \{quadratic\} & Y & Y \\
\hline
7 & \textbf{4, 5, 7, 11, 19} & 35 & 35 & Y & \{powers-of-2, constant\} & Y & N \\
8 & \textbf{1, 2, 10, 20, 100} & 200 & 200 & Y & \{powers-of-10\}, interleaved sequences & N & N \\
9 & \textbf{-2, 5, -4, 3, -6} & 1 & 1 & Y & \{linear, alternating\} & Y & N \\
10 & \textbf{1, 4, 9, 16, 25} & 36 & 36 & Y & \{quadratic\} & Y & Y \\
11 & \textbf{75, 15, 25, 5, 15} & 3 & 8.333 & N & \mysplit{\{reciprocal-powers-of-3, constant, alternating-linear\}\\$<$author: alternate(div by 5, add 10)$>$} & N & N \\
12 & \textbf{1, 2, 6, 24, 120} & 720 & - & N & \mysplit{-\\ $<$author: mult by 2,3,4,5,...$>$} & N & N \\
\hline
13 & \textbf{183, 305, 527, 749, 961} & 183 & 1162 & N & \mysplit{\{reciprocal-powers-of-10, quadratic\}\\ $<$VISUAL: each digit adds 2, mod 10$>$} & N & N \\
14 & \textbf{16, 22, 34, 58, 106} & 202 & 202 & Y & \{powers-of-2, constant\} & Y & N \\
15 & \textbf{17, 40, 61, 80, 97} & 112 & 112 & Y & \{quadratic\} & Y & Y \\
16 & \textbf{55, 34, 21, 13, 8} & 5 & 5 & Y & \{each term is the difference between the previous 2\} & N & N \\
17 & \textbf{259, 131, 67, 35, 19} & 11 & 11 & Y & \{reciprocal-powers-of-2, constant\} & N & N \\
18 & \textbf{93, 74, 57, 42, 29} & 18 & 18 & Y & \{quadratic\} & Y & Y \\
\hline
19 & \textbf{7, 21, 14, 42, 28} & 84 & 84 & Y & \{powers-of-2\}, interleaved sequences & N & N \\
20 & \textbf{2, -12, -32, -58, -90} & -128 & -128 & Y & \{quadratic\} & Y & Y \\
21 & \textbf{0, 9, 36, 81, 144} & 225 & 225 & Y & \{quadratic\} & Y & Y \\
22 & \textbf{15, 29, 56, 108, 208} & 401 & 401 & Y & \{each term is the sum of previous 4\} & N & N \\
23 & \textbf{13, -21, 34, -55, 89} & -144 & -144 & Y & \{each term is the difference between the previous 2\} & N & N \\
24 & \textbf{52, 56, 48, 64, 32} & 96 & -1024 & N & \mysplit{\{powers-of-2, each term is the difference between the previous 2, powers-of-8\}\\$<$author: alternate(add/subtract powers-of-2)$>$} & N & N \\
\hline

\end{tabular}
  \caption{Sequence problems from www.fibonicci.com, part I.
Successive columns show the index number of the problem, the sequence problem itself, the correct answer (per source), our answer (the DyMVeC answer), whether DyMVeC was correct, the explanation, whether example vectors \textbf{EMV1} and \textbf{EMV2} got the correct answer.  The Explanation entry in curly braces \{\} is the explanation given by  DyMVeC, if any.  In cases where  DyMVeC was wrong, following in angle brackets $<>$ is the interpretation of the sequence by the author.}
  \label{tab:fibprob1}
\end{table}

\begin{table}[hbp]
\hspace*{-3.5cm}
\begin{tabular}{|l|p{4.5cm}|l|l|l|p{7.5cm}|l|l|}
\hline
	& Problem & Ans & Ours & OK & Explanation & EMV1 & EMV2 \\
\hline
25 & \textbf{1, 3, 7, 11, 13} & 17 & 11 & N & \mysplit{\{cubic\}\\$<$author: cycle(add 2,4,4)$>$} & N & N \\
26 & \textbf{230, 460, 46, 92, 9.2} & 18.4 & 18.4 & Y & \{reciprocal-powers-of-5\}, interleaved sequences & N & N \\
27 & \textbf{3, 12, 24, 33, 66} & 75 & 121.5 & N & \mysplit{\{reciprocal-powers-of-2, each term is the sum of previous 3\}\\$<$author: alternate(add 9, mult by 2)$>$} & N & N \\
28 & \textbf{3, 5, 7, 11, 13} & 17 & 17 & Y & \{primes\} & Y & Y \\
29 & \textbf{8, 16, 24, 36, 48} & 64 & 64 & Y & \{quadratic, alternating\} & Y & N \\
30 & \textbf{68, 36, 20, 12, 8} & 6 & 6 & Y & \{reciprocal-powers-of-2, constant\} & N & N \\
\hline31 & \textbf{18, 6, 24, 8, 32} & 10.67 & 18 & N & \mysplit{\{each term is the sum of previous 2, constant, alternating\}\\$<$author: alternate(div by 3, mult by 4)$>$} & N & N \\
32 & \textbf{99, 18, 36, 9, 18} & 9 & -18 & N & \mysplit{\{each term is the difference between the previous 2, alternating-linear\}\\$<$VISUAL: alternate(add digits, mult by 2)$>$} & N & N \\
33 & \textbf{3, 8, 23, 68, 203} & 608 & 608 & Y & \{constant, powers-of-3\} & N & N \\
34 & \textbf{144, 73, 14, 8, 236} & 119 & - & N & \mysplit{-\\$<$author: cycle(add 2 and halve, then random)$>$} & N & N \\
35 & \textbf{10, 45, 15, 38, 20} & 31 & 31 & Y & \{linear, alternating-linear\} & N & N \\
36 & \textbf{1, 10, 37, 82, 145} & 226 & 226 & Y & \{quadratic\} & Y & Y \\
\hline
37 & \textbf{-2, 1, 6, 13, 22} & 33 & 33 & Y & \{quadratic\} & Y & Y \\
38 & \textbf{34, -21, 13, -8, 5} & -3 & -3 & Y & \{each term is the sum of previous 2\} & N & Y \\
39 & \textbf{1, 0, 1, -1, 2} & -3 & -3 & Y & \{each term is the difference between the previous 2\} & N & N \\
40 & \textbf{108, 56, 29, 15, 8} & 4 & 5 & N & \mysplit{\{reciprocal-powers-of-2, quadratic\}\\$<$author: add (4,2,1,1,0) then halve$>$} & N & N \\
41 & \textbf{1, 8, 27, 64, 125} & 216 & 216 & Y & \{cubic\} & N & N \\
42 & \textbf{-3, 3, 27, 69, 129} & 207 & 207 & Y & \{quadratic\} & Y & Y \\
\hline
\end{tabular}
  \caption{Sequence problems from www.fibonicci.com, part II.
Successive columns show the index number of the problem, the sequence problem itself, the correct answer (per source), our answer (the DyMVeC answer), whether DyMVeC was correct, the explanation, whether example vectors \textbf{EMV1} and \textbf{EMV2} got the correct answer.  The Explanation entry in curly braces \{\} is the explanation given by  DyMVeC, if any.  In cases where  DyMVeC was wrong, following in angle brackets $<>$ is the interpretation of the sequence by the author.}
  \label{tab:fibprob2}
\end{table}

\begin{table}[hbp]
\hspace*{-3.5cm}
\begin{tabular}{|l|p{4.5cm}|l|l|l|p{7.5cm}|l|l|}
\hline
	& Problem & Ans & Ours & OK & Explanation & EMV1 & EMV2 \\
\hline
43 & \textbf{3, 8, 15, 24, 35} & 48 & 48 & Y & \{quadratic\} & Y & Y \\
44 & \textbf{6, 14, 18, 28, 30} & 42 & 42 & Y & \{linear, alternating-linear\} & N & N \\
45 & \textbf{1/6, 1/3, 0.5, 2/3, 5/6} & 1 & 1 & Y & \{linear\} & Y & Y \\
46 & \textbf{4, 1, 0, 1, 4} & 9 & 9 & Y & \{quadratic\} & Y & Y \\
47 & \textbf{2, 8, 26, 80, 242} & 728 & 728 & Y & \{constant, powers-of-3\} & N & N \\
48 & \textbf{15, 210, 115, 220, 125} & 230 & 221 & N & \mysplit{\{reciprocal-powers-of-10, constant\}, interleaved sequences\\$<$VISUAL: alternate(1,2) prepended to multiples of 5$>$} & N & N \\
\hline
49 & \textbf{2, 7, 4, 14, 6} & 21 & 21 & Y & \{linear, alternating-linear\} & N & N \\
50 & \textbf{-1, 4, 1, 6, 3} & 8 & 8 & Y & \{linear, alternating\} & Y & N \\
51 & \textbf{1, 5, 13, 29, 61} & 125 & 125 & Y & \{powers-of-2, constant\} & Y & N \\
52 & \textbf{10, 21, 33, 46, 60} & 75 & 75 & Y & \{quadratic\} & Y & Y \\
53 & \textbf{88, 168, 248, 328, 408} & 488 & 488 & Y & \{linear\} & Y & Y \\
54 & \textbf{27, 214, 221, 228, 235} & 242 & - & N & \mysplit{-\\$<$VISUAL: 2 prepended to multiples of 7$>$} & N & N \\
\hline
55 & \textbf{0.78, 0.28, -0.22, -0.72, -1.22} & -1.72 & -1.72 & Y & \{linear\} & Y & Y \\
56 & \textbf{951, 1951, 950, 1950, 949} & 1949 & 1949 & Y & \{linear, alternating\} & Y & N \\
57 & \textbf{3, 5, 15, 17, 51} & 53 & 77 & N & \mysplit{\{powers-of-2, constant, alternating-linear\}\\$<$author: alternate(add 2, mult by 3)$>$} & N & N \\
58 & \textbf{15, 225, 16, 256, 17} & 289 & 287 & N & \mysplit{\{linear, alternating-linear\}\\$<$author: interleave(linear, quadratic)$>$} & N & N \\
59 & \textbf{66, 512, 618, 524, 630} & 536 & - & N & \mysplit{-\\$<$VISUAL: alternate(5,6) prepended to multiples of 6)$>$} & N & N \\
60 & \textbf{-100, 100, -200, 0, -300} & -100 & -100 & Y & \{linear, alternating\} & Y & N \\
\hline
61 & \textbf{3.3, 3.9, 3.27, 3.81, 3.243)} & 3.729 & - & N & \mysplit{-\\$<$VISUAL: 3. prepended to powers-of-3$>$} & N & N \\
62 & \textbf{100, 20, 2, 2/15, 1/150} & 1/3750 & - & N & \mysplit{-\\$<$author: div by multiples of 5$>$} & N & N \\
\hline
& Total Correct & & & 45 & & 28 & 19 \\
\hline
\end{tabular}
  \caption{Sequence problems from www.queendom.com. 
Successive columns show the index number of the problem, the sequence problem itself, the correct answer (per source), our answer (the DyMVeC answer), whether DyMVeC was correct, the explanation, whether example vectors \textbf{EMV1} and \textbf{EMV2} got the correct answer.  The Explanation entry in curly braces \{\} is the explanation given by DyMVeC, if any.  In cases where  DyMVeC was wrong, following in angle brackets $<>$ is the interpretation of the sequence by the author.}
  \label{tab:queenprob}
\end{table}

\hspace*{1cm}

\end{appendices}

 

\begin{thebibliography}{}
\bibitem{Strang} Gilbert Strang \textit{"Introduction to Linear Algebra, 5th Edition."}  Wellesley-Cambridge Press, 2016. 

 \end{thebibliography}
\end{document}